\documentclass{elsart3-1}

\usepackage{amssymb}

\usepackage[english,francais]{babel}

\newtheorem{theorem}{Theorem}[section]

\newtheorem{e-proposition}[theorem]{Proposition}

\newtheorem{e-definition}[theorem]{Definition\rm}

\renewcommand{\theequation}{\arabic{equation}}
\def\og{\leavevmode\raise.3ex\hbox{$\scriptscriptstyle\langle\!\langle$~}}
\def\fg{\leavevmode\raise.3ex\hbox{~$\!\scriptscriptstyle\,\rangle\!\rangle$}}

\def\wtil{\widetilde}

\newcommand{\om}{\omega}
\def\Pmu{\P_{\!\!\mu}}

\newcommand{\Sig}{\Sigma}

\newcommand{\R}{{\mathbb R}}

\def\N{{\mathbb N}}

\newcommand{\Prob}{{\mathbb P}\,}
\def\P{\Prob}
\newcommand{\Nat}{{\mathbb N}}

\def\be{\begin{equation}}
\def\ee{\end{equation}}

\def\ov{\overline}

\journal{the Acad\'emie des sciences}
\begin{document}
\centerline{}
\begin{frontmatter}


\selectlanguage{english}
\title{Hausdorff dimension of the multiplicative golden mean shift}


\selectlanguage{english}
\author[authorlabel1]{Richard Kenyon},
\ead{rkenyon@math.brown.edu}
\author[authorlabel2]{Yuval Peres},
\ead{peres@microsoft.com}
\author[authorlabel3]{Boris Solomyak}
\ead{solomyak@math.washington.edu}

\address[authorlabel1]{Richard Kenyon, Brown University}
\address[authorlabel2]{Yuval Peres, Microsoft Research}
\address[authorlabel3]{Boris Solomyak, University of Washington}


\medskip
\begin{center}
{\small Received *****; accepted after revision +++++\\
Presented by £££££}
\end{center}

\selectlanguage{english}
We compute the Hausdorff dimension of the ``multiplicative golden mean shift'' defined as the set of
all reals in $[0,1]$ whose binary expansion $(x_k)$ satisfies $x_k x_{2k}=0$ for all $k\ge 1$, and show that it is smaller than the Minkowski dimension.

\vskip 0.5\baselineskip

\selectlanguage{francais}
\noindent{\bf R\'esum\'e} \vskip 0.5\baselineskip \noindent
{\bf Dimension de Hausdorff du shift de Fibonacci multiplicatif. }
Nous calculons la dimension de Hausdorff
du ``shift de Fibonacci multiplicatif", c'est-\`a-dire l'ensemble
des nombres r\'eels dans $[0,1]$ dont le d\'eveloppement en binaire
$(x_k)$ satisfait $x_kx_{2k}=0$ pour tout $k\ge 1$. Nous
montrons que la dimension de Hausdorff est plus petite que la dimension de Minkowski.

\end{frontmatter}

\selectlanguage{francais}

\selectlanguage{english}

\section{Introduction}

A classical result of Furstenberg \cite{Furst} says that if $X$ is a closed subset of $[0,1]$, invariant under the
map $T_m:\ x\mapsto mx$ (mod 1), then its Hausdorff dimension equals the Minkowski (box-counting) dimension,
which equals the topological entropy of $T_m|_X$ divided by $\log m$. A simple example is the set
$
\Psi_{G}:= \Bigl\{ x = \sum_{k=1}^\infty x_k 2^{-k}:\ x_k \in \{0,1\},\ x_k x_{k+1}=0 \ \mbox{for all}\ k\Bigr\}
$
for which we have
$
\dim_H(\Psi_{G}) = \dim_M(\Psi_{G}) = \log_2\Bigl(\frac{1+\sqrt{5}}{2}\Bigr)
$
(the subscript $G$ here stands for the ``Golden Ratio'').
Instead, we consider the set
$$
\Xi_{G}:= \Bigl\{ x = \sum_{k=1}^\infty x_k 2^{-k}:\ x_k \in \{0,1\},\ x_k x_{2k}=0 \ \mbox{for all}\ k\Bigr\}
$$
which we call the ``multiplicative golden mean shift.'' The reason for this term is that the
set of binary sequences corresponding
to the points of $\Xi_{G}$ is invariant under the action of the semigroup of multiplicative positive integers $\N^*$:
$
M_r(x_k) = (x_{rk})\ \ \mbox{for}\ r\in \N.
$
Fan, Liao, and Ma \cite{Fan} showed that
$
\dim_M(\Xi_{G}) = \sum_{k=1}^\infty 
2^{-k-1}\log_2 F_{k+1}= 0.82429\ldots,
$
where $F_k$ is the $k$-th Fibonacci number: $F_1=1,\ F_2 = 2, F_{k+1} = F_{k-1}+F_k$, and raised the question of computing
the Hausdorff dimension of $\Xi_{G}$.

\medskip

\begin{theorem} \label{prop-gold}
We have  $\dim_H(\Xi_{G}) < \dim_M(\Xi_{G})$. In fact,
\be \label{gold1}
\dim_H(\Xi_{G}) = -\log_2 p = 0.81137\ldots,\ \ \mbox{where}\ p^3=(1-p)^2,\ \ \ 0<p<1.
\ee
\end{theorem}

Our manuscript \cite{KPS} contains substantial generalizations of this result, extending it to
 a large class of
``multiplicative subshifts.'' We state one of them at the end of the paper.

Although the set $\Xi_{G}$ is on the real line, it appears to have a strong resemblance with a class of self-affine sets
on the plane, namely, the Bedford-McMullen ``carpets'' \cite{Bedf,McM}, for which also the Hausdorff dimension is typically
smaller than the Minkowski dimension. However, this seems to be more of an analogy than a direct link.

An additional motivation to study the multiplicative subshifts comes from questions on multifractal analysis of multiple ergodic
averages raised in \cite{Fan}. Perhaps, the simplest non-trivial case of such multifractal analysis is the study of the sets
$
A_\theta:= \Bigl\{x = \sum_{k=1}^\infty x_k 2^{-k}:\ x_k \in \{0,1\},
\ \lim_{n\to \infty} \frac{1}{n} \sum_{k=1}^n x_k x_{2k} = \theta\Bigr\}\,.
$
It is not hard to show that $\dim_H(A_0) = \dim_H(\Xi_{G})$, which we compute in Theorem~\ref{prop-gold}. 
With more work, our method can be used to compute the Hausdorff dimension of $A_\theta$, but the details are beyond the scope of this note.

In this paper, we focus on $\Xi_{G}$ to explain our ideas and methods in the simplest possible setting.
To conclude the introduction, we should mention that the dimensions of some analogous sets, e.g.,
$
\wtil{\Xi} = \Bigl\{x= \sum_{k=1}^\infty x_k 2^{-k}:\ x_k \in \{0,1\},\ x_{k} x_{2k} x_{3k} = 0\ \mbox{for all $k$}\ \Bigr\}
$
are so far out of reach.

\section{Proof of Theorem~\ref{prop-gold}}

It is more convenient to work in the symbolic space $\Sig_2 = \{0,1\}^\N$, with the metric
$
\varrho((x_k) ,(y_k)) = \\ 2^{-\min\{n:\ x_n \ne y_n\}}.
$
It is well-known that the dimensions of a compact subset of $[0,1]$ and the corresponding set of binary digit sequences in $\Sig_2$
are equal (this is equivalent to replacing the covers by arbitrary interval with
those by dyadic intervals). Thus, it suffices to
determine the dimensions of the set $X_{G}$---the collection of all binary sequences $(x_k)$ such
that $x_k x_{2k}=0$ for all $k$.
Observe that
\be \label{eq1}
X_{G} = \Bigl\{\om =
{(x_k)}_{k=1}^\infty \in \Sig_2:\ {(x_{i2^r})}_{r=0}^\infty \in \Sig_{G}\ \ \mbox{for all}\ i \ \mbox{odd}\Bigr\}
\ee
where $\Sig_{G}$ is usual (additive) golden mean shift:
$
\Sig_{G}:= \{{(x_k)}_{k=1}^\infty\in \Sig_2,\ x_k x_{k+1}=0,\ \forall\,k\ge 1\}.
$

We will use the following well-known result; it essentially goes back to Billingsley \cite{Billing}.
We state it in the symbolic space $\Sig_2$ where $[u]$ denotes the cylinder set of sequences starting with
a finite ``word'' $u$ and $x_1^n = x_1\ldots x_n$.

\medskip

\begin{prop}[see \cite{Falc}] \label{prop-mass}
Let $E$ be a Borel set in $\Sig_2$ and let $\nu$ be a finite Borel measure on $\Sig_2$.

{\bf (i)} If
$\liminf_{n\to \infty} (-\frac{1}{n}) \log_2 \nu[x_1^n] \ge s\ \ \mbox{for $\nu$-a.e.}\ x\in E,$
then $\dim_H(E) \ge s$.

{\bf (ii)} If
$\liminf_{n\to \infty} (-\frac{1}{n}) \log_2 \nu[x_1^n] \le s\ \ \mbox{for all}\ x\in E,$
then $\dim_H(E) \le s$.
\end{prop}

\medskip

Given a probability measure $\mu$ on $\Sig_{G}$, we can define a probability measure on $X_{G}$ by
\be \label{eq-meas1}
\Pmu[u]:= \prod_{i\le n,\, i\ \mbox{\tiny odd}} \mu[u|_{J(i)}], \ \ \mbox{where}\ J(i) = \{2^r i\}_{r=0}^\infty
\ee
and $u|_{J(i)}$ denotes the ``restriction'' of the word $u$ to the subsequence $J(i)$.
It turns out that this class of measures is sufficiently rich to compute $\dim_H(X_{G})$.

For $k\ge 1$ let
$\alpha_k$ be the partition of $\Sig_{G}$ into cylinders of length $k$.
For a measure $\mu$ on $\Sig_2$ and a finite partition $\alpha$, denote by $H^\mu(\alpha)$ the $\mu$-entropy of the
partition, with base $2$ logarithms:
$
H^\mu(\alpha) = -\sum_{C\in \alpha} \mu(C)\log_2\mu(C).
$
Define
\be\label{def-smu}
s(\mu):= \sum_{k=1}^\infty \frac{H^\mu(\alpha_k)}{2^{k+1}}\,.
\ee

\begin{prop} \label{prop-ldim}
Let $\mu$ be a probability measure on $\Sig_{G}$. Then $\dim_H(X_{G})\ge s(\mu)$.
\end{prop}

\medskip

\noindent{\bf Proof.}
We are going to demonstrate that for every $\ell \in \N$,
\be \label{eq-lb2}
\liminf_{n\to \infty} \frac{-\log_2\Pmu[x_1^n]}{n} \ge
\sum_{k=1}^\ell \frac{H^\mu(\alpha_k)}{2^{k+1}}\ \ \mbox{for $\Pmu$-a.e.}\ x.
\ee
Then, letting $\ell\to\infty$ and using Proposition~\ref{prop-mass}(i) will yield
the desired inequality.
Fix $\ell\in \Nat$. By a routine argument, to verify (\ref{eq-lb2}) we can restrict ourselves to $n=2^\ell r, \ r\in \N$.
In view of (\ref{eq-meas1}), we have
\be \label{eq-lb3}
\Pmu[x_1^n] \le \prod_{k=1}^\ell \ \prod_{\frac{n}{2^k} < i \le \frac{n}{2^{k-1}},\ i\ \mbox{\tiny odd}}
\mu[x_1^n|_{J(i)}].
\ee
Note that $x_1^n|_{J(i)}$ is a word of length $k$ for $i\in (n/2^k, n/2^{k-1}]$, with $i$ odd, which is
a beginning of a sequence in $\Sig_G$.
Thus, $[x_1^n|_{J(i)}]$ is an element of the partition $\alpha_k$. 
The random variables $x\mapsto -\log_2\mu[x_1^n|_{J(i)}]$ are i.i.d\
for $i\in (n/2^k, n/2^{k-1}]$, with $i$ odd,
and their expectation  equals $H^\mu(\alpha_k)$, by the definition of entropy.
Note that there are $n/2^{k+1}$ odds in $(n/2^k, n/2^{k-1}]$.
Fixing $k,\ell$ with $k\le \ell$ and taking $n=2^\ell r$,
$r\to \infty$, we get an infinite sequence of i.i.d.\ random variables. Therefore, by a version of the Law of Large
Numbers,

\be \label{eq-lb4}
\forall\ k\le \ell, \sum_{\frac{n}{2^k} < i \le \frac{n}{2^{k-1}},\ i\ \mbox{\tiny odd}}
\frac{-\log_2 \mu[x_1^n|_{J(i)}]}{(n/2^{k+1})} \
\to H^\mu(\alpha_k)\ \ \mbox{as}\ n = 2^\ell r\to \infty,\ \ \mbox{for $\Pmu$-a.e.\ $x$}.
\ee
By (\ref{eq-lb3}) and (\ref{eq-lb4}), for $\Pmu$-a.e.\ $x$,
$$
\frac{-\log_2\Pmu[x_1^n]}{n} \ge  \sum_{k=1}^\ell \frac{1}{2^{k+1}} 
\sum_{\frac{n}{2^k} < i \le \frac{n}{2^{k-1}},\ i\ \mbox{\tiny odd}}
\frac{-\log_2\mu[x_1^n|_{J(i)}]}{n/2^{k+1}} \to \sum_{k=1}^\ell \frac{H^\mu(\alpha_k)}{2^{k+1}}\,.
$$
This confirms (\ref{eq-lb2}), so the proof is complete.
\qed

\medskip

\noindent{\bf Proof of the lower bound for the Hausdorff dimension in Theorem~\ref{prop-gold}.}
Let $s:= \sup\{s(\mu):\ \mu$ is a probability measure on $\Sig_G\}$.
By Proposition~\ref{prop-ldim}, we have $\dim_H(X_G) \ge s$, and we will prove that this is actually an equality.
To this end, we specify a measure which will turn out to be ``optimal.'' This measure is Markov,
but non-stationary. It could be ``guessed'' or derived by solving the optimization problem (which also yields that the optimal measure is unique). However, for the proof of dimension formula it suffices to produce the answer.
Let $\mu$ be a Markov measure on $\Sig_G$, with initial probabilities $(p,1-p)$, and the matrix of transition probabilities
$P = (P(i,j))_{i,j=0,1} = \left( \begin{array}{cc} p & \ 1-p \\ 1 & 0 \end{array} \right)$. 
Using elementary properties of entropy, it is
not hard to see that $s(\mu) = \frac{H(p)}{2} + \frac{ps(\mu)}{2} + \frac{(1-p)s(\mu)}{4}$, whence $s(\mu) = \frac{2H(p)}{3-p}$. Maximizing over $p$  yields $s(\mu) = 2\log_2\frac{p}{1-p}$, and comparing this to $s(\mu) = \frac{2H(p)}{3-p}$ we get
\be \label{eq-s2}
p^3=(1-p)^2,\ \ s(\mu)=-\log_2 p.
\ee
Combined with Proposition~\ref{prop-ldim}, this proves the lower bound for the Hausdorff dimension in
(\ref{gold1}). \qed

\medskip

\noindent {\bf Proof of the upper bound for the Hausdorff dimension in Theorem~\ref{prop-gold}.}
Denote by $N_i(u)$ the number of symbols
$i$ in a word $u$. By the definition of the measure $\mu$,  
we obtain for any $u=u_1\ldots u_k \in \{0,1\}^n$,
\be
\mu[u]   =  p_{u_1} P(u_1,u_2)\cdot \ldots \cdot P(u_{k-1}, u_k)
         =  (1-p)^{N_1(u_1\ldots u_k)} p^{N_0(u_1\ldots u_k) - N_1(u_1\ldots u_{k-1})}. \label{eq-meas3}
\ee
Indeed, the probability of a 1 is always $1-p$,  whereas the probability of a 0 is $p$, except in those cases when it follows a 1, which it must
by the definition of $\Sig_G$.
In view of (\ref{eq-meas3}), by the definition of the measure $\Pmu$ on $X_G$, we have 
$
\Pmu[x_1^n] = (1-p)^{N_1(x_1^n)} p^{N_0(x_1^n) - N_1(x_1^{n/2})}
$ for any $x\in X_G$ and $n$ even.
Using that $(1-p)^2=p^3$ and $N_0(x_1^n) = n - N_1(x_1^n)$, we obtain that
$$
\Pmu[x_1^n] = p^n p^{N_1(x_1^n)/2- N_1(x_1^{n/2})}.
$$
Let
$a_\ell= -\frac{1}{n}\log_2 \Pmu[x_1^n]$ for $n=2^\ell$. 
Then
$
a_\ell = -\log_2 p\left( 1+ \frac{1}{2} \Bigl[\frac{N_1(x_1^n)}{n} - \frac{N_1(x_1^{n/2})}{n/2}\Bigr]\right).
$
Now we see that the average of $a_\ell$'s ``telescopes'':
$$
\frac{a_1 + \cdots + a_\ell}{\ell} = -\log_2 p \left(1+ \frac{1}{2\ell} \Bigl[\frac{N_1(x_1^{2^\ell})}{2^\ell} - N_1(x_1)\Bigr]\right)
\to -\log_2 p, \ \ \mbox{as}\ \ell\to \infty.
$$
It follows that
$$
\liminf_{\ell\to \infty} a_\ell = \liminf_{\ell\to \infty} 2^{-\ell} (-\log_2 \Pmu[x_1^{2^\ell}]) \le -\log_2 p = s,
$$
for every $x\in X_G$, so $\dim_H(X_G) \le s$ by Proposition~\ref{prop-mass}(ii). \qed

\section{Generalization}

Here we state a generalization of Theorem~\ref{prop-gold} to the case of arbitrary multiplicative subshifts of
finite type; the proof can be found in \cite{KPS}.

\medskip

\begin{theorem} \label{th-main} {\rm (i)} Let $A$ be a $0$-$1$ primitive $m\times m$ matrix
(i.e.\ some power of $A$ has only positive entries).
Consider
$
\Xi_A = \Bigl\{ x = \sum_{k=1}^\infty x_k m^{-k}:\ x_k \in \{0,\ldots,m-1\},\ A(x_k, x_{2k})=1 \ \mbox{for all}\ k\Bigr\}.
$
Then
$
\dim_H(\Xi_A) = \frac{1}{2}\log_m \sum_{i=0}^{m-1} t_i,
$
where $(t_i)_{i=0}^{m-1}$ is the unique vector satisfying
$
t_i^2= \sum_{j=0}^{m-1} A(i,j) t_j,\ \ t_i>1,\ i=0,\ldots,m-1.
$

{\rm (ii)} The Minkowski dimension of $\Xi_A$ exists and equals
$
\dim_M(\Xi_A) = \sum_{k=1}^\infty 2^{-k-1} \log_m(A^{k-1} \ov{1},\ov{1})$ where $\ov{1}=(1,\ldots,1)^T \in \R^m.
$
We have $\dim_H(\Xi_A) = \dim_M(\Xi_A)$ if and only if all row sums of $A$ are equal.
\end{theorem}





\section*{Acknowledgements}

We are grateful to J\"org Schmeling for telling us about the problem, and to Aihua Fan, Lingmin Liao, and Jihua Ma for
sending us their preprint \cite{Fan} prior to publication.
The research of R. K. and
B. S.  was supported in part by NSF.

\end{document}